\newcommand{\fancyplain}{\plain}
\newcommand{\mathbb}[1]{{\bf #1}}

\documentclass[fancyheadings,twoside,psfig,fleqn,epsf,12pt]{article}
\usepackage [dvips]{graphicx}
\usepackage{psfig}
\usepackage{epsf}
\usepackage{amssymb}
\usepackage{fancyheadings}

\newcommand{\be}{\begin{equation}}
\newcommand{\ee}{\end{equation}}
\newcommand{\bea}{\begin{eqnarray}}
\newcommand{\eea}{\end{eqnarray}}
\newcommand{\bean}{\begin{eqnarray*}}
\newcommand{\eean}{\end{eqnarray*}}

\newcommand{\bp}{\backprime}
\newcommand{\pp}{\prime}
\newcommand{\OPT}{\mbox{\tiny OPT}}
\newcommand{\LL}{\mbox{\tiny L}}
\newcommand{\RR}{\mbox{\tiny R}}
\newcommand{\CC}{\mbox{\tiny C}}
\newcommand{\NE}{\mbox{\tiny NE}}
\newcommand{\NW}{\mbox{\tiny NW}}
\newcommand{\SE}{\mbox{\tiny SE}}
\newcommand{\SW}{\mbox{\tiny SW}}
\newcommand{\un}{\bar{u}^{n}}
\newcommand{\uj}{\bar{u}_{j}^{n}}
\newcommand{\ujp}{\bar{u}_{j+1}^{n}}
\newcommand{\ujm}{\bar{u}_{j-1}^{n}}

\newenvironment{remarks}{{\flushleft \bf Remarks:}}{}

\newenvironment{acknowledgment}{{\flushleft \bf Acknowledgment:}}{}

\begin{document}

\title{Compact Central WENO Schemes\\for Multidimensional Conservation Laws}
\author{Doron Levy\footnotemark[2] \and Gabriella Puppo\footnotemark[3] \and 
Giovanni Russo\footnotemark[4]}
\date{}

\renewcommand{\thefootnote}{\fnsymbol{footnote}}
\footnotetext[2]
{Department of Mathematics, University of California, Berkeley, CA 94720, and Lawrence
Berkeley National Lab; dlevy@math.berkeley.edu}
\footnotetext[3]{Dipartimento di Matematica, Politecnico di Torino, 
 Corso Duca degli Abruzzi 24, 10129 Torino, Italy; puppo@calvino.polito.it}
\footnotetext[4]{Dipartimento di Matematica, Universit\`{a} dell'Aquila, 
 Via Vetoio, loc. Coppito - 67100 L'Aquila, Italy; russo@univaq.it}
\renewcommand{\thefootnote}{\arabic{footnote}}

\maketitle


\begin{abstract}
We present a new third-order central scheme 
for approximating solutions of systems of conservation laws
in one and two space dimensions.  In the spirit of Godunov-type schemes, 
our method is based on reconstructing a piecewise-polynomial interpolant
from cell-averages which is then advanced exactly in time.  

In the reconstruction step, we introduce a new third-order, {\em compact},
CWENO reconstruction, which is written
as a convex combination of interpolants based on different stencils.
The heart of the matter is that 
one of these interpolants is taken as an arbitrary 
quadratic polynomial and the weights of the convex combination
are set as to obtain third-order accuracy in smooth regions.  
The embedded mechanism in the WENO-like schemes guarantees
that in regions with discontinuities or large gradients, there
is an automatic switch to a one-sided second-order reconstruction,
which prevents the creation of spurious oscillations.

In the one-dimensional case, our new third order scheme is based on
an extremely compact {\em four\/} point stencil.  Analogous
compactness is retained in more space dimensions.  The
accuracy, robustness and high-resolution properties of 
our scheme are demonstrated in a variety of one and
two dimensional problems.

\end{abstract}


\bigskip
\noindent
{\bf AMS(MOS) subject classification.} Primary 65M10; secondary 65M05.



\section{Introduction}			\label{section:introduction}
\setcounter{equation}{0}
\setcounter{figure}{0}
\setcounter{table}{0}

We are concerned with multidimensional
systems of hyperbolic conservation laws of the form
\be
  u_{t} + \nabla_{x} \cdot f(u) = 0,  \qquad x \in \mathbb{R}^{d}, \quad
u=(u_{1},\ldots,u_{n}).
 \label{eq:conservation}
\ee

Methods for approximating solutions to equation (\ref{eq:conservation})
have attracted a lot of attention in recent years 
(see \cite{god-rav:difference}, \cite{leveque:numerical-methods}, 
\cite{tadmor:approximate} and the references therein).

In this work we focus on Godunov-type schemes, where
one first reconstructs a piecewise-polynomial 
interpolant which is then advanced exactly in time according
to (\ref{eq:conservation}) and finally projected on its cell-averages.   
Generally, one can divide Godunov-type schemes into two sub-classes -
{\em upwind methods\/} and {\em central methods}.

In {\em upwind schemes}, one first reconstructs a polynomial in every cell,
which is then used to compute a new cell average in the same location in the 
next time step.  This procedure requires solving Riemann
problems at the discontinuous interfaces.  For high-order
methods, instead of analytically solving the resulting Riemann
problems, one typically implements approximate Riemann solvers
or some form of flux splitting.
For systems of conservation laws, or in the demanding context
of more space dimensions, this procedure turns out to be more intricate
as such Riemann solvers do not exist. 
The essentially non-oscillatory (ENO) methods by Harten are 
the prototype of high-order methods (see \cite{harten-eoc:eno},
 \cite{shu-osher:eno} and the references therein). A recent review 
of ENO and WENO methods can
be found in \cite{shu:eno}.

{\em Central schemes}, on the other hand, are based on averaging over
the Riemann fans, a procedure which is typically done by staggering
between two grids.  They require no Riemann solvers, 
no projection along characteristic directions and no flux splitting.
Therefore, all that one has to do in order to solve a problem is to supply the
flux function.  They are more simple when compared with upwind schemes.

The major difference between
different central methods is in the reconstruction step, where one computes a
piecewise-polynomial interpolant from the previously computed cell-averages. 

The prototype of central schemes is the Lax-Friedrichs \cite{lxf}
scheme which is based on a piecewise-constant interpolant. 
Even though it is very robust, it is only first-order accurate
and, moreover, it suffers from excessive numerical dissipation.  
A second-order central method was proposed by Nessyahu and
Tadmor in \cite{nessyahu-tadmor:non-oscillatory}.  
This method is based on a MUSCL-like \cite{van-leer:ultimate}
piecewise linear interpolant and nonlinear limiters which
prevent spurious oscillations (see \cite{sanders-weiser} for a
different approach).  
A variety of extensions to these methods were suggested.
The one-dimensional third-order method of Liu and Tadmor 
\cite{liu-tadmor:3rd} is  based on the  third-order reconstruction by
Liu and Osher in \cite{liu-osher:nonosc}.  
For the two-dimensional method see 
\cite{arminjon:2d} and \cite{jiang-tadmor:nonosc}.

A first step for importing the high-order reconstructions that were
derived in the upwind framework was taken
in \cite{bianco-puppo-russo:central}.  There, the ENO method was
transformed into the central setup, and a new {\em mostly centered stencil\/}
was shown to produce the least oscillatory results.  The next
step was taken in the 1D case in \cite{levy-puppo-russo:1d}, where
a new {\em central weighted non-oscillatory\/} (CWENO) reconstruction was
introduced. This CWENO method is based on the upwind WENO
methods by \cite{liu-osher-chan:weno} and \cite{jiang-shu:weno}, in which 
an interpolant is written as a convex combination of several
reconstructions which are based on different stencils.
These methods include an internal switch which is designed such as
to provide the maximum possible accuracy in smooth regions, while 
automatically switching to a the more robust one-sided stencil in 
the presence of discontinuities and large-gradients.
The 1D CWENO method was extended
to two space dimensions in \cite{levy-puppo-russo:2d-3} and 
\cite{levy-puppo-russo:2d-4}. For a  numerical study of the behavior
of the total variation for the CWENO scheme, see
\cite{levy-puppo-russo:tv}.

In this paper we present a new, compact CWENO reconstruction.  This new
reconstruction is based on defining an arbitrary quadratic function which
is added to linear interpolants in such a way as to obtain third-order
accuracy in smooth regions (in one and two space dimensions).
In regions with discontinuities or large gradients, the weights
are automatically changed so that they switch to a one-sided
second-order linear reconstruction.
These reconstructions turn to be extremely compact; in the one
dimensional case, e.g., the reconstruction is based on a four-point stencil.

The structure of this paper is as follows:  
We start in \S\ref{section:balance} with a brief overview of central schemes for
conservation laws in one and two space dimensions.  

We then proceed  to present our new, compact, third-order
CWENO reconstruction in \S\ref{section:cweno}.  
The idea of introducing an arbitrary quadratic reconstruction
is first presented in the one dimensional framework and then extended to two space
dimensions.

We end in \S\ref{section:examples} where we demonstrate our new method in 
several test cases.   First, the accuracy tests (both in one-dimensional 
and two-dimensional cases) show the third-order accuracy of the method.
We then solve the one-dimensional system of the Euler equations of 
gas dynamics for a few test problems 
and we illustrate the behavior of the scheme in
scalar two-dimensional cases.  In particular, we would
like to stress that in our numerical results we observe a very robust
and non-oscillatory behavior of the weights, which can be related to the
overall robustness and accuracy properties of our new method.

\begin{acknowledgment}
The work of D.L. was supported in part by the Applied Mathematical Sciences 
subprogram of the Office of Energy Research, U.S. Department of Energy, under 
contract DE--AC03--76--SF00098.
Part of this work was done while G.P. and G.R. were visiting the Lawrence Berkeley Lab.
\end{acknowledgment}


\section{Central Schemes for Conservation Laws}		\label{section:balance}
\setcounter{equation}{0}
\setcounter{figure}{0}
\setcounter{table}{0}

In this section we give a brief overview of central schemes for approximating 
solutions to hyperbolic conservation laws in one and two space dimensions.  
For further details we refer the reader to 
\cite{tadmor:approximate}, \cite{jiang-tadmor:nonosc}, 
\cite{levy-puppo-russo:1d} and the references therein.

Starting in the one-dimensional case, we seek numerical solutions of the
Cauchy problem
\be
\left\{
\begin{array}{l}
 \displaystyle{u_{t} + f(u)_{x} = 0,}   \\ \\
 \displaystyle{u(x,t\!=\!0) = u_{0}(x).} \label{eq:1d.conservation}
\end{array} \right.
\ee
For simplicity, we introduce a uniformly spaced grid in the $(x,t)$ space, where
the mesh spacings are denoted by $h:=\Delta x$ and $k:=\Delta t$, respectively.
We denote by $\bar{u}_{j}^{n}$, the numerical approximation of the cell average
in the cell $I_{j} := [x_{j-1/2},x_{j+1/2}]$ at time $t^{n} = nk$,
where $x_{j} = j h$.  The finite-difference method will approximate the
cell-averages at time $t^{n+1}$ based on their values at time $t^{n}$.

We start by reconstructing at time $t^{n}$ a piecewise-polynomial 
conservative interpolant from the known cell-averages, $\bar{u}_{j}^{n}$, i.e.,
\be
 P_{u}(x,t^{n}) := \sum_{j} R_{j}(x) \chi_{j},    \label{eq:1d.reconst}
\ee
where $\chi_{j}$ is the characteristic function of the interval $I_{j}$
and $R_{j}(x)$ is a polynomial defined in $I_{j}$.  It is
here, in the reconstruction step, where the accuracy and non-oscillatory requirements 
enter.  Different central methods will be typically based on different reconstructions.

The reconstruction, $P_{u}(x,t^{n})$, is then evolved exactly in time (integrating
(\ref{eq:1d.conservation})), and then projected on staggered cells, in order to compute
the cell average at $I_{j+1/2}$.  With this procedure we obtain
\be
 \bar{u}_{j+1/2}^{n+1} =  \bar{u}_{j+1/2}^{n}
 + \frac{1}{h} \int_{t^{n}}^{t^{n+1}}[f(P_u(x_{j},\tau)) - f(P_u(x_{j+1},\tau))]\, d\tau. 
  \label{eq:1d.exact} 
\ee
Based on the reconstruction, (\ref{eq:1d.reconst}), the first term
on the RHS of (\ref{eq:1d.exact}) can be explicitly computed,
\[
 \bar{u}_{j+1/2}^{n} = \frac{1}{h} \int_{x_{j}}^{x_{j+1}}P_{u}(x,t^{n})\, dx.
\]
In order to compute the second integral on the RHS of (\ref{eq:1d.exact}),
namely the integral in time over the fluxes, one should observe that
due to the staggering, up to a suitable CFL condition, these integrals
involve only smooth functions, and hence  can be approximated by a sufficiently
smooth quadrature.  A second-order method can be obtained, e.g., using
the mid-point rule in time (see \cite{nessyahu-tadmor:non-oscillatory}). 
Simpson's rule for the quadrature in time will provide 
fourth order accuracy (which will naturally be sufficient also for
a third-order method).

The quadratures for the time integrals of 
the fluxes, require the {\em prediction\/} of point-values
of the function in several points in the interval $[t^{n},t^{n+1}]$.  One
possible approach is to use a Taylor expansion based on the
equation, (\ref{eq:1d.conservation}).  Such an approach was used, e.g., 
in \cite{nessyahu-tadmor:non-oscillatory} and \cite{liu-tadmor:3rd}.  
In order to avoid the technical complications involved in the Taylor
series (in particular with high-order methods, and when dealing with systems
of equations), one can alternatively use a Runge-Kutta (RK) method 
directly on (\ref{eq:1d.conservation}), to predict the required values
in later times.  
By using the {\em Natural Continuous Extension\/} (NCE) of Runge-Kutta schemes
\cite{zennaro:rk}, 
the value of the flux at the nodes of the quadrature 
formula can be computed with a single RK step. 
Such a method is simpler compared with the Taylor
expansion method, but it does require another reconstruction: the 
reconstruction of the point values of the derivatives of the fluxes
at time $t^{n}$ which are then used as the input to the RK solver
(see \cite{bianco-puppo-russo:central} and \cite{levy-puppo-russo:1d} for more
details).

The simplicity of central schemes manifests itself when turning
to deal with systems of equations.  Basically, the algorithm that
was described in the scalar case, repeats itself component-wise.
Based on the type of the reconstruction, when solving systems
of equations, there can even be
simplifications over a purely componentwise extension
of the scalar scheme.  These can be found, e.g., in 
\S\ref{subsection:systems}
below.

The extension to two space dimensions is straightforward;
We consider
\be
 u_{t} + f(u)_{x} + g(u)_{y} = 0,   \label{2d.conser}
\ee
subject to the initial condition, $u(x,y,t\!=\!0) = u_{0}(x,y)$.
Here, $\Delta x$ and $\Delta y$ will denote the spatial mesh spacings, 
while $\Delta t$ denotes the spacing in time.  The two-dimensional
cells are now $I_{i,j} = [x_{i-1/2},x_{i+1/2}] \times [y_{j-1/2},y_{j+1/2}]$.

Following the general methodology, we wish to construct the cell-averages,
 $\bar{u}_{j,k}^{n+1}$, based on the cell averages at time $t^{n}$.  First,
we construct a two-dimensional interpolant which reads
\be
 P_{u}(x,y,t^{n}) = \sum_{i,j} R_{i,j}(x,y) \chi_{i,j},  \label{eq:2d.interp}
\ee
with $\chi_{i,j}$ being the characteristic function of the cell $I_{i,j}$
and $R_{i,j}(x,y)$ a polynomial of a suitable degree.  An exact evolution
in time of the interpolant (\ref{eq:2d.interp}) which is 
projected on its cell-averages now reads (compare with (\ref{eq:1d.exact}))
\bea
 \bar{u}_{i+1/2,j+1/2}^{n+1} &=&  \bar{u}_{i+1/2,j+1/2}^{n} +   \label{eq:2d.exact.evolution} \\
&& + \frac{1}{\Delta x \Delta y} \int_{t^{n}}^{t^{n+1}} \int_{y=y_{j}}^{y_{j+1}}
 [f(P_u(x_{i},y,\tau)) - f(P_u(x_{i+1},y,\tau)]\, dy d\tau + \nonumber \\
&& + \frac{1}{\Delta x \Delta y} \int_{t^{n}}^{t^{n+1}} \int_{x=x_{i}}^{x_{i+1}}
 [f(P_u(x,y_{j},\tau)) - f(P_u(x,y_{j+1},\tau)] \, dx d\tau \nonumber.
\eea
The staggered cell-average at time $t^{n}$ can be directly computed by
\[
 \bar{u}_{i+1/2,j+1/2}^{n} = \frac{1}{\Delta x\Delta y} 
 \int_{x_{i}}^{x_{i+1}} \int_{y_{j}}^{y_{j+1}} P_{u}(x,y,t^{n})\, dy dx.
\]
Analogously to the one-dimensional case, the flux integrals on the 
RHS of (\ref{eq:2d.exact.evolution}) can be approximated using a quadrature 
rule in time.  This can be coupled with a quadrature rule for the line 
integrals in space. Here it is necessary to avoid quadrature points at which 
the fluxes may not be smooth.
The details can be found, e.g., in \cite{levy-puppo-russo:2d-3,levy-puppo-russo:2d-4}
(see also \cite{jiang-tadmor:nonosc}).


\section{A Compact third-order CWENO Reconstruction}	\label{section:cweno}
\setcounter{equation}{0}
\setcounter{figure}{0}
\setcounter{table}{0}


\subsection{The One-Dimensional Framework}	\label{subsection:1d}

In this section we derive our new CWENO reconstruction in one space dimension.
The two dimensional extension will follow in \S\ref{subsection:2d}.

We first note that in the absence of large gradients, we obtain third order
 accuracy if we choose for the reconstruction the {\em optimal} polynomial
\[
 P_j(x) = P_{\OPT ,j}(x),
\]
where $P_{\OPT ,j}(x)$ is the parabola that interpolates the data
 $\ujm , \uj ,\ujp$ in the sense of cell averages to enforce conservation:
\[ \int_{x_{j+l-1/2}}^{x_{j+l+1/2}} P_{\OPT ,j}(x)\, dx = 
      \bar{u}_{j+l}^{n}, \qquad l = -1,0,1.
\]
These conditions determine $P_{\OPT ,j}(x)$ completely, namely:
\be
 P_{\OPT ,j}(x) = u_{j}^{n} + u_{j}^{\pp}(x-x_{j}) 
 + \frac{1}{2} u_{j}^{\pp\pp}(x-x_{j})^{2},   \label{eq:1d.p.exact}
\ee
with
\bean
 u_{j}^{n} & = & \uj - \frac{1}{24}(\ujp-2\uj + \ujm), \\
 \qquad u_{j}^{\pp} & = & \frac{\ujp - \ujm}{2 \Delta x}, \qquad
 u_{j}^{\pp \pp}   =  \frac{\ujm - 2\uj + \ujp}{\Delta x^2}.
\eean

However, when discontinuities or large gradients occur, this
reconstruction would be oscillatory.
Therefore following the WENO methodology (\cite{liu-osher-chan:weno}, 
\cite{jiang-shu:weno}, \cite{levy-puppo-russo:1d}), 
we construct an essentially non-oscillatory
interpolant as a convex combination of polynomials which are
based on different stencils.  
Specifically, in the cell $I_{j}$ we write
\be
 P_{j}(x) = \sum_{i} w_{i}^{j} P_{i}^{j}(x),
  \quad \sum_{i}w^j_{i} = 1, \quad w_{i} \geq 0, \quad i \in \{\LL,\CC,\RR\},
 \label{eq:1d.combination}
\ee
where $P_{\LL}$ and $P_{\RR}$ are linear functions based on a left stencil 
and a right
stencil, respectively, and $P_{\CC}$ is a quadratic polynomial.
In order to simplify the notations, we will omit the
upper index $j$, remembering that the weights and the three polynomials change
from cell to cell.

Conservation requires that $P_{\RR}(x)$ will interpolate the cell averages
\[
 \int_{x_{j-1/2}}^{x_{j+1/2}} P_{\RR}(x) dx = \Delta x \, \uj, \qquad
 \int_{x_{j+1/2}}^{x_{j+3/2}} P_{\RR}(x) dx = \Delta x \, \ujp,
\]
which in turn, results with
\be
 P_{\RR}(x) = \uj + \frac{\ujp - \uj}{\Delta x}(x-x_{j}).  \label{eq:1d.p.r}
\ee
Similarly, for the left interpolant we have
\be
 P_{\LL}(x) = \uj + \frac{\uj - \ujm}{\Delta x}(x-x_{j}).  \label{eq:1d.p.l}
\ee
All that is left is to reconstruct a centered polynomial, $P_{\CC}$,
such that the convex combination, (\ref{eq:1d.combination}),
will be third-order accurate in smooth regions.  It must, therefore, satisfy
\be
 P_{\OPT}(x) = C_{\LL} P_{\LL}(x) + C_{\RR} P_{\RR}(x) 
  + C_{\CC} P_{\CC}(x),   \quad
 \sum_{i} C_{i} = 1, \quad i \in \{\LL,\CC,\RR\},
\label{eq:exact.combination}
\ee
where $C_{\LL}, C_{\CC}$ and $C_{\RR}$ are constants.
Due to the staggering between every two consecutive steps of the central
method, our reconstruction should provide half-cell averages which are
third-order accurate.  A straightforward calculation shows that {\em any\/}
symmetric choice of constants $C_{i}$ in (\ref{eq:exact.combination})
provides the desired accuracy.
In particular, for the specific choice of $C_{\LL} \!=\! C_{\RR} \!=\! 1/4$, 
equations (\ref{eq:1d.p.r})--(\ref{eq:exact.combination}) yield
\bea
 P_{\CC}(x) & = & 2 P_{\OPT}(x) - \frac{1}{2} (P_{\RR}(x) + P_{\LL}(x)) =
        \uj - \frac{1}{12}(\ujp - 2\uj + \ujm) + \nonumber \\
   && + \frac{\ujp-\ujm}{2 \Delta x}(x-x_{j}) + 
        \frac{\ujp-2\uj+\ujm}{\Delta x^2}(x-x_j)^{2}.   \label{eq:1d.center}
\eea

In order to complete the reconstruction of $P_{j}(x)$ in
(\ref{eq:1d.combination}), it is left to compute the weights $w_{i}$.
Following \cite{jiang-shu:weno}, \cite{levy-puppo-russo:1d}, we write
\be
 w_{i} = \frac{\alpha_{i}}{\sum_{k}\alpha_{k}}, 
 \qquad \alpha_{i} = \frac{C_{i}}{(\varepsilon + I\!S_{i})^{p}},
 \qquad i,k \in \{\LL,\CC,\RR\}.  \label{eq:weights}
\ee
The constants $C_{i}$'s in (\ref{eq:weights}) are chosen to be the
same as in (\ref{eq:exact.combination}), i.e.,
$C_{\LL} = C_{\RR}= 1/4$, while  $C_{\CC} = 1/2$.

The {\em smoothness indicators\/}, $I\!S_{i}$, are 
responsible for detecting large gradients or discontinuities and to automatically
switch to the stencil that generates the least oscillatory reconstruction in
such cases. 
Once again, we follow \cite{jiang-shu:weno}, \cite{levy-puppo-russo:1d} and define
in each cell $I_{j}$ the three smoothness indicators, $I\!S_{i}$, as
\be
 I\!S_{i} = \sum_{l=1}^{2} \int_{x_{j-1/2}}^{x_{j+1/2}} h^{2l-1}(P_{i}^{(l)}(x))^2 dx.
   \qquad i \in \{\LL,\CC,\RR\}. \label{eq:1d.smoothness.indicators}
\ee
A direct computation,
based on (\ref{eq:1d.p.r}), (\ref{eq:1d.p.l}) and (\ref{eq:1d.center}), yields
\bea
 I\!S_{\LL} & = & (\uj - \ujm)^2, \qquad  I\!S_{\RR} = (\ujp - \uj)^2,  \nonumber \\
 I\!S_{\CC} & = & \frac{13}{3}(\ujp-2\uj+\ujm)^2 + \frac{1}{4}(\ujp - \ujm)^2.
  \label{eq:si.1d}
\eea
The constant $\varepsilon$ is taken as to prevent the denominator from vanishing.
Furthermore, its valued has to satisfy two requirements, namely 
\begin{description}
\item{i)} $\varepsilon \gg IS$ in smooth regions
\item{ii)} $\varepsilon \ll IS$ near discontinuities
\end{description}
The first conditions guarantees that, on smooth regions, the weights are basically equal to
the constants that provide high accuracy. 
The second conditions guarantees that   
in the presence of a discontinuity, the weights of the parabola and of
one of the one-sided linear reconstructions will be practically zero, 
and we will be left with the other one-sided linear reconstruction.
Hence, our third-order method automatically switches to a second-order
method in the presence of large gradients, which is exactly what
makes it so robust as will be evident in our numerical computations presented below.

The constant $p$ weights the departure from smoothness. 
We used the value $p=2$. For a discussion about its choice see 
\cite{jiang-shu:weno} and \cite{levy-puppo-russo:1d}.

A second, non-oscillatory reconstruction is required for the flux
derivative.  It is natural to adapt the  reconstruction that we used
for the half-cell averages also to the reconstruction of the
point-values of the flux derivative.  
Here however the interpolation requirements will be in the sense of
point values instead of cell averages.
Once again, a direct computation shows
that any symmetric choice of constants will provide the desired accuracy and
it will only be natural to use the same constants, $C_{i}$'s, that
were used in (\ref{eq:exact.combination}).
This is simpler than the case of \cite{levy-puppo-russo:1d} where
we had to use different sets of constants for the two different reconstructions
(the reconstruction of the half cell averages, and the reconstruction of
the point-values of the flux derivative at the edges of the domain).

\begin{remarks}
\begin{enumerate}
\item 
We would like to emphasize that our new method is based on adding the
arbitrary parabola $P_{\CC}$ into the convex combination which is the heart
of our non-oscillatory reconstruction.   Our numerical simulations
showed that the freedom we have in selecting the constants $C_{i}$ 
has no influence on the  properties of the method.  
It is easy to prove that we obtain
third-order accuracy regardless of the smoothness of the weights, 
as long as they are symmetric.  This is substantially more robust
than the third-order method of Liu and Tadmor in \cite{liu-tadmor:3rd}, where
the order of the method did depend on the smoothness of the limiters,
and could deteriorate to first order in the presence of large gradients
(as was shown in \cite{bianco-puppo-russo:central}).
\item By extending our new ideas, one can modify our previous
CWENO method presented in \cite{levy-puppo-russo:1d} in order to obtain
a fifth-order, central, non-oscillatory scheme.  This can be done
by simply adding a fourth-order polynomial for the computation of
the point-values.
\item Our new reconstruction is equivalent to limiting with the minimum
slope instead of slope zero in the presence of a discontinuity.
Hence, it is based on a second-order reconstruction close to shocks, unlike 
the scheme of \cite{liu-tadmor:3rd} which can be only first order accurate
in such regions.
\item In the one dimensional case, the additional parabola is needed
only for the accurate recovery of the point-values.  In the 2D framework, 
however, the equivalent additional parabola will be required also for
the accurate reconstruction of the fractional cell-averages.
\end{enumerate}
\end{remarks}


\subsection{A Two-Dimensional Extension} 		\label{subsection:2d}

We extend the ideas of \S\ref{subsection:1d} to the two-dimensional framework.
This extension is straightforward and is based on reconstructing an interpolant
as a convex combination of four one-sided, piecewise-linear interpolants, and
a centered,  quadratic interpolant, such as to get the desired third-order
accuracy in smooth regions.  

\begin{figure}[htbp]    
\begin{center}
  \leavevmode {
  \hbox{
        \epsfxsize=5cm
        \epsffile{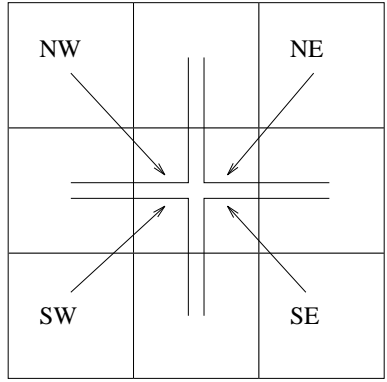}
       }
  }
\end{center}
\caption{{\sf The Two-Dimensional Stencil}\label{figure:2d_stencil}}
\end{figure}

Following these ideas, the reconstruction in the cell ${I_{ij}}$, can be written as
\be
 P_{i,j}(x,y) = \sum_{k} w_{k}^{i,j} P_{k}^{i,j}(x,y),
  \qquad k \in \{\NE, \NW, \SE, \SW, \CC \},
 \label{eq:2d.combination}
\ee
with $\sum_{k} w_{k}^{i,j} = 1$, and $w_{k}^{i,j} \geq 0$.
Here, $P_{\NE}$, $P_{\NW}$, $P_{\SE}$ and $P_{\SW}$ 
are the one-sided linear reconstructions, and $P_{\CC}$ is
a centered quadratic reconstruction (see Figure \ref{figure:2d_stencil}).
Similar to the one-dimensional case, we will simplify our notations by
omitting the superscripts, remembering that both the weights and
the polynomials, change from cell to cell.

Clearly, in an analog to the one-dimensional case, 
(\ref{eq:1d.p.r}--\ref{eq:1d.p.l}),  the four linear reconstructions are given by
\bea
P_{\NE}(x,y) & = & \un_{i,j} + \frac{\un_{i+1,j}-\un_{i,j}}{\Delta x}(x-x_i) 
	+ \frac{\un_{i,j+1}-\un_{i,j}}{\Delta y}(y-y_j),  \nonumber \\
P_{\NW}(x,y) & = & \un_{i,j} + \frac{\un_{i,j}-\un_{i-1,j}}{\Delta x}(x-x_i) 
	+ \frac{\un_{i,j+1}-\un_{i,j}}{\Delta y}(y-y_j),  \nonumber \\
P_{\SW}(x,y) & = & \un_{i,j} + \frac{\un_{i,j}-\un_{i-1,j}}{\Delta x}(x-x_i) 
	+ \frac{\un_{i,j}-\un_{i,j-1}}{\Delta y}(y-y_j),  \nonumber \\
P_{\SE}(x,y) & = & \un_{i,j} + \frac{\un_{i+1,j}-\un_{i,j}}{\Delta x}(x-x_i) 
	+ \frac{\un_{i,j}-\un_{i,j-1}}{\Delta y}(y-y_j).
 \label{eq:2d.linear.reconstructions}
\eea
The centered polynomial, $P_{\CC}(x,y)$, is taken such as to satisfy
\be
 P_{\OPT}(x,y) = \sum_{k} C_{k} P_{k}(x,y),
  \quad \sum_{k}C_{k}=1,
  \qquad k \in \{\NE, \NW, \SE, \SW, \CC \}.   \label{eq:2d.comb}
\ee
Here, $P_{\OPT}$ is the quadratic polynomial based on a nine-point stencil,
centered around $I_{i,j}$, which is given by (see \cite{levy:third})
\bea
 P_{\OPT} (x,y) & = &
  \tilde{u}_{i,j}^{n} + \tilde{u}_{i,j}^{\pp} (x-x_i) + \tilde{u}_{i,j}^{\bp}(y-y_j) +
  \tilde{u}_{i,j}^{\pp \bp} (x-x_i)(y-y_j) + \nonumber \\
 && +\frac{1}{2}\tilde{u}_{i,j}^{\pp \pp}(x-x_i)^2 + \frac{1}{2}\tilde{u}_{i,j}^{\bp \bp}(y-y_j)^2,
 \label{eq:2d.exact}
\eea
where
\bea
 \tilde{u}_{i,j}^{n} &=& \bar{u}_{i,j} - \frac{1}{24}\left( (\Delta x)^2 \tilde{u}_{i,j}^{\pp \pp} 
  + (\Delta y)^2 \tilde{u}_{i,j}^{\bp \bp} \right),   \nonumber \\
 \tilde{u}_{i,j}^{\pp} &=& \frac{\un_{i+1,j}-\un_{i-1,j}}{2 \Delta x}, \qquad
 \tilde{u}_{i,j}^{\bp} = \frac{\un_{i,j+1}-\un_{i,j-1}}{2 \Delta y}, \nonumber \\
 \tilde{u}_{i,j}^{\pp \pp} &=& \frac{\un_{i+1,j}-2\un_{i,j}+\un_{i-1,j}}{\Delta x^2}, \qquad
 \tilde{u}_{i,j}^{\bp \bp} = \frac{\un_{i,j+1}-2\un_{i,j}+\un_{i,j-1}}{\Delta y^2}, \\
 \tilde{u}_{i,j}^{\pp \bp} &=& \frac{\un_{i+1,j+1}+\un_{i-1,j-1}-\un_{i+1,j-1}-\un_{i-1,j+1}}{4 \Delta x \Delta y}. \nonumber
\eea
Unlike the one-dimensional case, not every symmetric selection of the
constants $C_{k}$'s will provide a third-order reconstruction for the
quarter cell-averages.  Here, a straightforward computation shows that 
in order to satisfy the accuracy requirements, we must take
 $C_{\NE} = C_{\NW} = C_{\SW} = C_{\SE} = 1/8$.  Hence, $C_{\CC}=1/2$, and 
(\ref{eq:2d.comb}) implies
\bea
 P_{\CC}(x,y) &=& 2 P_{\OPT}(x,y) - \frac{1}{4} 
  \bigg[ P_{\NE}(x,y) + P_{\NW}(x,y) + P_{\SW}(x,y) + P_{\SE}(x,y) \bigg] = \nonumber \\
 & = &  u_{i,j}^{n} + u_{i,j}^{\pp} (x-x_i) + u_{i,j}^{\bp}(y-y_j) +
  u_{i,j}^{\pp \bp} (x-x_i)(y-y_j) + \nonumber \\
 && +\frac{1}{2}u_{i,j}^{\pp \pp}(x-x_i)^2 + \frac{1}{2}u_{i,j}^{\bp \bp}(y-y_j)^2,
 \label{eq:2d.center}
\eea
where
\bean
 u_{i,j}^{n} &=& \bar{u}_{i,j}^{n}
  - \frac{1}{12} \bigg[ (\Delta x)^2 u_{i,j}^{\pp \pp} 
  + (\Delta y)^2 u_{i,j}^{\bp \bp} \bigg],   \\
 u_{i,j}^{\pp} &=& \tilde{u}_{i,j}^{\pp}, \quad
 u_{i,j}^{\bp} = \tilde{u}_{i,j}^{\bp}, \\
 u_{i,j}^{\pp \pp} &=& 2\tilde{u}_{i,j}^{\pp \pp}, \quad
 u_{i,j}^{\bp \bp} = 2\tilde{u}_{i,j}^{\bp \bp}, \quad
 u_{i,j}^{\pp \bp} = 2\tilde{u}_{i,j}^{\pp \bp}.
\eean

All that remains is to determine the weights $w_{k}^{i,j}$ in  (\ref{eq:2d.combination}).
Once again, we write
\[
 w_{k}^{i,j} = \frac{\alpha_{k}^{i,j}}{\sum_{l}\alpha_{l}^{i,j}}, 
 \qquad \alpha_{k}^{i,j} = \frac{C_{k}^{i,j}}{(\varepsilon + I\!S_{k}^{i,j})^{p}},
 \qquad k,l \in \{\NE,\NW,\SE,\SW,\CC\}.
\]
The constants, $C_{k}$, are the same constants that were used to reconstruct
the centered parabola in (\ref{eq:2d.comb}).
The constants, $\varepsilon$ and $p$, play the same role as in the one-dimensional
case.  
At that point, to simplify the notations we assume that the mesh spacings are 
equal in the $x$ and  $y$ directions, i.e., $\Delta x = \Delta y = h$.  We can then follow
\cite{levy-puppo-russo:2d-3}, and define 
the smoothness indicators, $I\!S_{k}^{i,j}$, as
\be
 I\!S^{i,j}_{k} = \sum_{|\alpha|=1,2} \int_{x_{i-h/2}}^{x_{i+h/2}}
  \int_{y_{j-h/2}}^{y_{j+h/2}} h^{2(|\alpha|-1)}(D^{\alpha}P_{k})^{2}, \qquad
  k \in \{\NE,\NW,\SE,\SW,\CC\}.
 \label{eq:2d.smoothness.indicators}
\ee
If $\Delta x \ne \Delta y$, then only a trivial enhancement to
(\ref{eq:2d.smoothness.indicators}) is required.
For the four one-sided linear reconstructions, which can be all written as
\[
 P_{k} = \hat{u} + \hat{u}^{\pp}(x-x_{i}) + \hat{u}^{\bp}(y-y_{j}),
\qquad k \in \{\NE,\NW,\SE,\SW\}.
\]
with suitable reconstructed point-values and first derivatives, 
a direct computation of (\ref{eq:2d.smoothness.indicators}) results with
\be
 I\!S_{k} = h^{2}[(\hat{u}^{\pp})^{2}+(\hat{u}^{\bp})^{2}].
\ee
The centered smoothness  indicator, $I\!S_{\CC}$, which corresponds to the 
centered quadratic reconstruction, $P_{c}(x,y)$, (\ref{eq:2d.center}), is given by
\[
 I\!S_{\CC} = h^{2} 
  \left[ (u^{\pp})^{2} + (u^{\bp})^{2} \right] 
   + \frac{h^{4}}{3} 
    \left[ 13(u^{\pp \pp})^{2} + 14 (u^{\pp \bp})^{2} + 13 (u^{\bp \bp})^{2}\right],
\]
(the discrete derivatives are given by (\ref{eq:2d.center})).

\subsection{A Note on Systems}   \label{subsection:systems}
Almost nothing changes when turning to deal with systems.  The reconstruction
that was described in the previous sections, directly transforms to systems and
is performed component-wise.  The only relatively subtle issue is the computation
of the smoothness indicators.

In our previous work \cite{levy-puppo-russo:1d}, we have suggested several approaches
for the computation of the smoothness indicators.  Three different options were
suggested: the first is to allow every component to have a strictly individual behavior, 
namely to allow a different stencil with different smoothness indicators to each 
component.  In the second approach, we designed {\em global smoothness indicators\/}
such as to force every component to adjust even when the discontinuity is in a 
different component.  The last approach was to use external information about
the system.  For example, in the Euler equations of gas dynamics, one expects 
both shocks and contact 
discontinuities to occur in the density.  Hence, all stencils can be 
tuned according to this component.  

Our results in \cite{levy-puppo-russo:1d} showed that the best approach is
the one which was based on the {\em global smoothness indicators}.  It requires
no additional information on the system, it produced less oscillatory results
compared with the individual smoothness indicators for each component, and
it was the simplest to implement.

In the one-dimensional case, e.g., the global smoothness indicators are given by
(compare with (\ref{eq:1d.smoothness.indicators}))
\be
 I\!S_k = \frac{1}{d} \sum_{r=1}^{d} \frac{1}{\|\bar{u}_r\|_2} 
       \left( \sum_{l=1}^{2}
       \int_{x_{j-1/2}}^{x_{j+1/2}} h^{2l-1} \left( P_{k,r}^{(l)}
            \right)^2 \; dx \right), \quad k\in \{\LL,\CC,\RR\}   \label{eq:IS.system}
\ee
Here the $k$-th polynomial in the $r$-th component is denoted by $P_{k,r}$,
and $d$ is the number of equations.
The scaling factor $\|\bar{u}_r\|_2$ is defined as the $L^2$ norm of the cell averages
of the $r$-th component of $u$, 
\[
 \|\bar{u}_r\|_2 = \left( \sum_{\mbox{\small all}\, j} |\bar{u}_{j,r}|^2 h \right)^{1/2}.
\]

The numerical examples performed
for systems, appearing in \S\ref{section:examples}, were carried out with 
 the global smoothness indicators.
	

\section{Examples}				\label{section:examples}
\setcounter{equation}{0}
\setcounter{figure}{0}
\setcounter{table}{0}


We present numerical tests in one and two space dimensions, in which we 
demonstrate the accuracy, robustness and high-resolution properties
of our new method.

Following our previous works (\cite{bianco-puppo-russo:central}, 
\cite{levy-puppo-russo:1d}, \cite{levy-puppo-russo:2d-3} and 
\cite{levy-puppo-russo:2d-4}), in all our numerical examples we 
integrate in time using a Runge-Kutta method with natural continuous 
extension, which was presented in \cite{zennaro:rk}, with a fixed time step. 

The time step is determined by imposing that the Courant number is a given 
fraction of the maximum Courant number determined by linear stability analysis.
The Courant number is defined by 
\[
   C = \frac{\lambda}{\max_j \rho_j}
\]
where $\rho_j$ denotes the spectral radius of the matrix $f'(u_j)$ computed on the 
initial condition, and $\lambda = \Delta t/\Delta x$ is the mesh ratio.

The linear stability analysis carried out in \cite{bianco-puppo-russo:central}
yields a Courant number $C = 3/7$ for the one
dimensional case. We remark that the stencil used in 
\cite{bianco-puppo-russo:central} was different than the one that we use, 
and therefore linear stability should be repeated for the compact scheme 
in order to obtain a sharp estimate on the maximum Courant number.

\subsection*{Example 1: Accuracy Tests}			\label{subsection:accuracy}

Our first example checks the accuracy of our new method in several one and
two-dimensional test cases.  In all of the one-dimensional tables, the norms of the
errors are given by
\[
   \begin{array}{ll}
      L^1-\mbox{error}: & 
            ||\mbox{Error}||_1 = \sum_{j=1}^N|u(x_j,t^n)-u_j^n| h, \\
      L^\infty-\mbox{error}: & 
            ||\mbox{Error}||_{\infty} = \max_{1\leq j\leq N} |u(x_j,t^n)-u_j^n|.
   \end{array}
\]
Analogous expressions hold for the two-dimensional norms.

\begin{enumerate}
 \item {\bf Linear advection:} This test estimates the convergence
rate at large times.   We solve $u_{t} + u_{x} = 0$, subject to the initial data
 $u(x,t\!=\!0)=\sin(\pi x)$ and to periodic boundary conditions on $[-1,1]$.
The integration time was taken as $T=10$.  
 
 In Table \ref{table:accuracy.sin.2}, we show the results obtained for this test
problem with $\epsilon=10^{-2}$.  We clearly see a third order convergence
rate in both $L^{1}$ and $L^{\infty}$ norms.  
We also show in Table \ref{table:accuracy.sin.6}
the results obtained for the same example with $\epsilon=10^{-6}$,
which is the value suggested in both \cite{jiang-shu:weno} and 
\cite{levy-puppo-russo:1d}.

Compared with Table \ref{table:accuracy.sin.2}, the errors are larger,
and the convergence rate is not as regular as before.  This is mainly due to the 
fact that for a very small value of $\varepsilon$, condition i) is not satisfied, 
until the grid spacing $\Delta x$ becomes very small. 

\item {\bf Linear advection with oscillatory data:}  
This test is used to detect deteriorations of accuracy due to
oscillations in the parameters that control the selection of the stencil
(for details see \cite{bianco-puppo-russo:central} and the references therein).
Once again, the equation
is $u_{t}+u_{x}=0$, subject to the oscillatory initial data, 
 $u(x,t\!=\!0)=\sin^{4}(\pi x)$ and to periodic boundary conditions on
 $[-1,1]$.  Here, the integration time is taken as $T=1$ and $\epsilon=10^{-2}$.

The results of this test are displayed in Table \ref{table:accuracy.sin4},
and confirm the third-order accuracy of the method with no deteriorations
in its accuracy.

\item {\bf Burgers equation:} We solve the Burgers equation
 $u_t + (0.5 u^2)_x = 0$, subject to the initial data
 $u(x,t\!=\!0) = 1+0.5\sin(\pi x)$ and to periodic boundary conditions on $[-1,1]$.
 The integration time is $T=.33$, and $\epsilon=10^{-2}$.  Here,
a shock develops at $T=2/\pi$. Note that here the maximum speed of
propagation is $f'(u) = 3/2$. Thus we use 
$\lambda = .66 * 3/7 \simeq 2/3\lambda_{\mbox{\tiny max}}$

Table \ref{table:accuracy.burgers} shows the results we obtained which
verify the third-order accuracy of the method also for nonlinear problems.

\item {\bf 2D Linear advection:} Finally, we implemented our method for
the two-dimensional linear advection problem, $u_t+u_x+u_y=0$.  The 
initial condition is taken as $u(x,t\!=\!0) = \sin^{2}(\pi x) \sin^{2}(\pi y)$,
and we impose periodic boundary conditions on $[0,1]^{2}$.  The errors
and the estimated convergence rate our computed at time $T=1$.

In Table \ref{table:accuracy.2d.constant} we present the results obtained
when the weights are taken as constants (\ref{eq:2d.comb}).
Table  \ref{table:accuracy.2d.nonconstant} shows the results obtained with
the fully non-linear scheme, where the weights of the reconstruction
include also the oscillatory indicators.  With constant weights our
method is third-order as expected, while with non-constant weights, there
seems to be a better convergence rate.  A careful study of the tables
shows, however, that this better convergence rate is mainly due to
larger errors on the coarser grids.

\end{enumerate}

\begin{table}[htbp]
\begin{center}
\begin{tabular}{||l||c|c|c|c||}
\hline \hline
N & $L^{1}$ error & $L^{1}$ order & $L^{\infty}$ error & $L^{\infty}$ order \\ \hline
   20 &   0.1423     &     -    & 0.1484     &   -       \\ 
   40 &   0.1308E-01 &   3.44 &   0.1708E-01 &   3.12 \\
   80 &   0.7054E-03 &   4.21 &   0.1071E-02 &   4.00 \\
  160 &   0.7517E-04 &   3.23 &   0.7823E-04 &   3.78 \\
  320 &   0.9391E-05 &   3.00 &   0.7977E-05 &   3.29 \\
  640 &   0.1174E-05 &   3.00 &   0.9406E-06 &   3.08 \\
 1280 &   0.1467E-06 &   3.00 &   0.1158E-06 &   3.02 \\ \hline \hline
\end{tabular}
\caption{Linear advection, $T=10$, $u_0(x) = \sin(\pi x)$, $\epsilon = 10^{-2}$,
Courant number $C=0.9 C_{\mbox{\tiny max}}$, $C_{\mbox{\tiny max}} = 3/7$
\label{table:accuracy.sin.2}}
\end{center}
\end{table}

\begin{table}[htbp]
\begin{center}
\begin{tabular}{||l||c|c|c|c||}
\hline \hline
N & $L^{1}$ error & $L^{1}$ order & $L^{\infty}$ error & $L^{\infty}$ order \\ \hline

   20 &   0.2292     &    -     &   0.2522     &    -      \\ 
   40 &   0.8975E-01 &   1.35 &   0.9943E-01 &   1.34 \\
   80 &   0.2184E-01 &   2.04 & 0.3759E-01 &   1.40 \\
  160 &   0.3677E-02 &   2.57 & 0.1090E-01 &   1.79 \\
  320 &   0.3682E-03 &   3.32 & 0.1896E-02 &   2.52 \\
  640 &   0.2454E-04 &   3.91 & 0.1585E-03 &   3.58 \\
 1280 &   0.1379E-05 &   4.15 & 0.5972E-05 &   4.73 \\ \hline \hline
\end{tabular}
\caption{Linear advection, $T=10$, $u_0(x) = \sin(\pi x)$, $\epsilon = 10^{-6}$,
Courant number $C=0.9 C_{\mbox{\tiny max}}$, $C_{\mbox{\tiny max}} = 3/7$
\label{table:accuracy.sin.6}}
\end{center}
\end{table}

\begin{table}[htbp]
\begin{center}
\begin{tabular}{||l||c|c|c|c||}
\hline \hline
N & $L^{1}$ error & $L^{1}$ order & $L^{\infty}$ error & $L^{\infty}$ order \\ \hline
   20 &   0.1285     &     -     &   0.1909     &     -    \\
   40 &   0.2813E-01 &   2.19  &   0.5223E-01 &   1.87 \\
   80 &   0.2608E-02 &   3.43  &   0.5226E-02 &   3.32 \\
  160 &   0.2553E-03 &   3.35  &   0.3619E-03 &   3.85 \\
  320 &   0.3055E-04 &   3.06  &   0.3319E-04 &   3.45 \\
  640 &   0.3826E-05 &   3.00  &   0.3814E-05 &   3.12 \\
 1280 &   0.4777E-06 &   3.00  &   0.4654E-06 &   3.03 \\ \hline \hline
\end{tabular}
\caption{Linear advection, $T=1$, $u_0(x) = \sin^{4}(\pi x)$, $\epsilon = 10^{-2}$, 
Courant number $C=0.9 C_{\mbox{\tiny max}}$, $C_{\mbox{\tiny max}} = 3/7$
\label{table:accuracy.sin4}}
\end{center}
\end{table}

\begin{table}[htbp]
\begin{center}
\begin{tabular}{||l||c|c|c|c||}
\hline \hline
N & $L^{1}$ error & $L^{1}$ order & $L^{\infty}$ error & $L^{\infty}$ order \\ \hline
   20 &   0.7974E-02 &      -   & 0.1527E-01 &        - \\
   40 &   0.6654E-03 &   3.58 &   0.1844E-02 &   3.05 \\
   80 &   0.6563E-04 &   3.34 &   0.2340E-03 &   2.98 \\
  160 &   0.8494E-05 &   2.95 &   0.3645E-04 &   2.68 \\
  320 &   0.1067E-05 &   2.99 &   0.4937E-05 &   2.88 \\
  640 &   0.1355E-06 &   2.98 &   0.6388E-06 &   2.95 \\
 1280 &   0.1695E-07 &   3.00 &   0.8047E-07 &   2.99 \\  \hline \hline
\end{tabular}
\caption{Burgers equation, $T=.33$, $u_0(x) =1 + 0.5 \sin(\pi x) $, $\epsilon = 10^{-2}$, 
$\lambda=0.66*3/7$, corresponding to a Counant number $C=0.99 C_{\mbox{\tiny max}}$
\label{table:accuracy.burgers}}
\end{center}
\end{table}

\begin{table}[htbp]
\begin{center}
\begin{tabular}{||l||c|c|c|c||}
\hline \hline
N & $L^{1}$ error & $L^{1}$ order & $L^{\infty}$ error & $L^{\infty}$ order \\ \hline
 10 &   3.696E-02  & -     & 1.252E-01 & - \\
 20 &   4.964E-03  &  2.90 & 1.767E-02 & 2.83 \\
 40 &   6.304E-04  &  2.98 & 2.264E-03 & 2.96 \\
 80 &   7.902E-05  &  3.00 & 2.842E-04 & 2.99 \\
 160 &  9.880E-06  &  3.00 & 3.555E-05 & 3.00 \\ 
  \hline \hline
\end{tabular}
\caption{2D Linear Advection, Constant weights; $T=1$, $\lambda = 0.425$, 
  $u_0(x) = \sin^{2}(\pi x)\sin^{2}(\pi y) $,
\label{table:accuracy.2d.constant}}
\end{center}
\end{table}

\begin{table}[htbp]
\begin{center}
\begin{tabular}{||l||c|c|c|c||}
\hline \hline
N & $L^{1}$ error & $L^{1}$ order & $L^{\infty}$ error & $L^{\infty}$ order \\ \hline
 10 &   9.750E-02  & -     & 3.447E-01 & - \\
 20 &   1.419E-02  &  2.78 & 8.111E-02 & 2.09 \\
 40 &   9.387E-04  &  3.92 & 7.967E-03 & 3.35 \\
 80 &   8.319E-05  &  3.50 & 4.465E-04 & 4.16 \\
 160 &  9.977E-06  &  3.06 & 3.999E-05 & 3.48 \\ 
  \hline \hline
\end{tabular}
\caption{2D Linear Advection, Non-Constant weights; $T=1$, $\lambda = 0.425$, 
  $u_0(x) = \sin^{2}(\pi x)\sin^{2}(\pi y) $, $\epsilon = 10^{-2}$
\label{table:accuracy.2d.nonconstant}}
\end{center}
\end{table}


\subsection*{Example 2: 1D Systems - Euler Equations of Gas Dynamics}

\begin{figure}[htbp]
  \begin{center}
  \includegraphics[width=4.5in]{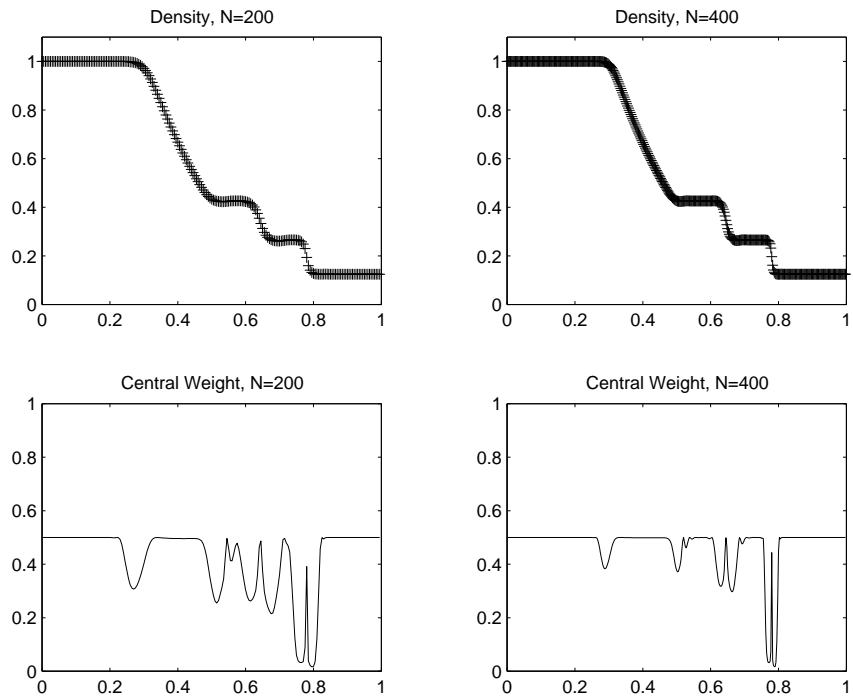}
  \end{center}
  \caption{Euler equations of gas dynamics - Sod initial data, $\lambda=0.1$, $T=0.16$
 \label{figure:sod}}
\end{figure}

\begin{figure}[htbp]
  \begin{center}
  \includegraphics[width=4.5in]{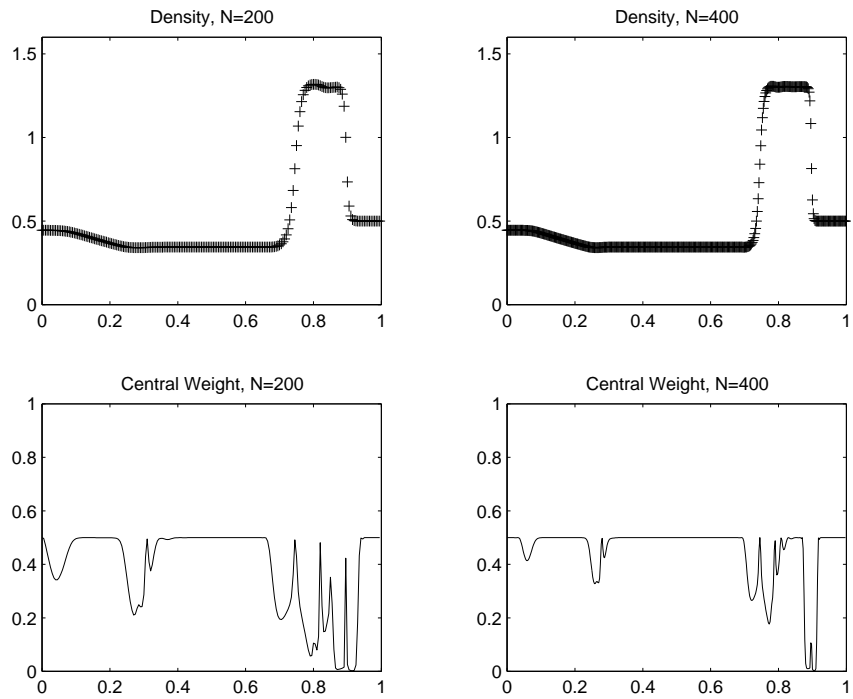}
  \end{center}
  \caption{Euler equations of gas dynamics - Lax initial data, $\lambda=0.1$, $T=0.16$
 \label{figure:lax}}
\end{figure}

Next, we solve the Euler equations of gas dynamics, 
\[
\frac{\partial}{\partial t} \left( \begin{array}{l} \rho\\ m\\ E \end{array}
\right)+
\frac{\partial}{\partial x} \left( \begin{array}{c} m\\ \rho u^2+p\\ u(E+p)
\end{array} \right)=0,
\qquad p=(\gamma-1) \cdot \left(E-\frac{\rho}{2}u^2\right),
\]
where the variables
$\,\rho,\, u,\, m=\rho u,\, p\,$ and $\,E\,$ denote the density, velocity,
momentum, pressure and total energy, respectively. 
We use two sets of initial data:

\begin{enumerate}
\item Shock tube problem with Sod's initial data, \cite{sod:hyperbolic},
\[
  \left\{  \begin{array}{lr}
    (\rho_l, m_l, E_l) = (1,0,2.5), &  x < 0.5, \\  
    (\rho_r, m_r, E_r) = (0.125,0,0.25),  &  x > 0.5.
    \end{array}  \right.
\]

\item Shock tube problem with Lax initial data, \cite{lax:weak},
\[
  \left\{  \begin{array}{lr}
    (\rho_l, m_l, E_l) = (0.445,0.311,8.928), &  x < 0.5, \\  
    (\rho_r, m_r, E_r) = (0.5,0,1.4275),  &  x > 0.5.
    \end{array}  \right.
\]
\end{enumerate}
We integrate the equations up to time $T=0.16$.
In Figure \ref{figure:sod} we show the density components for Sod's
initial data, and in Figure \ref{figure:lax} we show the equivalent
plot for Lax initial data.  In both figures we also present the weight of
the central stencil at the final time.  All the results are given for 200 and 400
grid points. Following \cite{nessyahu-tadmor:non-oscillatory}, we pick 
$\lambda = .1$.
Note that the maximum characteristic speed for Sod's problem is roughly $2.5$,
while for Lax problem the maximum propagation speed is $\simeq 5$. 
Since our scheme has a Courant number no larger than $.5$, we see that
$\lambda = .1$ is actually the maximum value compatible with stability.
This explains while there are still some wiggles in the test
solution of Lax, while the Sod's solution is monotone.

It is interesting to compare the behavior of the central weight of the new method
to the behavior of the central weight of the original CWENO method \cite{levy-puppo-russo:1d}.
Here, the weights are much smoother compared with the weights in \cite{levy-puppo-russo:1d}.
The accuracy and stability properties of the method can be related to
the smoothness of the nonlinear weights involved (see \cite{bianco-puppo-russo:central}).

\subsection*{Example 3: 2D Problems}
\begin{enumerate}
\item {\bf Linear rotation:}
Following \cite{levy-puppo-russo:2d-3}, we consider a linear rotation
of a square patch on $[0,1]^{2}$, with initial condition
 $u_{0}(x,y) = 1$ for $\{|x-1/2|\leq 1/2\} \times \{|y-1/2|\leq 1/2\}$ and
zero elsewhere.  In Figure \ref{figure:rotation} we display the
solution after a rotation of $\pi/4$ and of $\pi/2$.  There are no 
spurious oscillations.  We also show the corresponding
central weight.  As expected, this weight is zero in regions where 
the solution has steep gradients, and that is exactly the property
that prevents spurious oscillations from developing.  
Even though we are dealing with linear waves, the
resulting resolution is relatively good.  Due to the compactness of 
the stencil, when the slopes are not sharp, they
are not identified as discontinuities. This can be observed in 
the plots of the central weight, which returns to its constant value (1/2)
on the slope.

\begin{figure}[htbp]
  \begin{center}
  \includegraphics[width=4.5in]{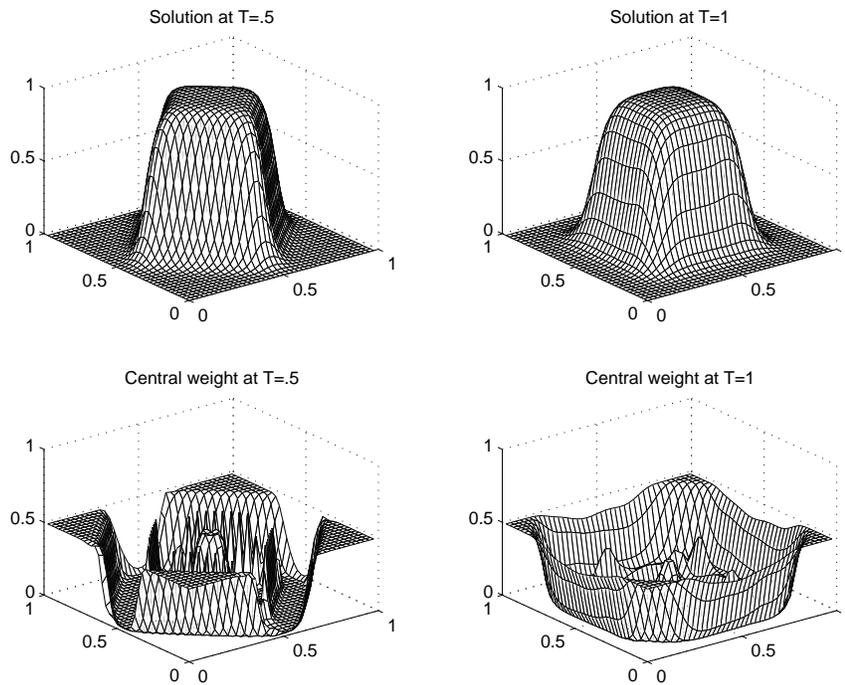}
  \end{center}
  \caption{Linear Rotation, $\lambda=0.425$, $N=40$
 \label{figure:rotation}}
\end{figure}

\item {\bf 2D - Burgers equation:}
We end by solving the two-dimensional Burgers equation, 
 $u_t + (u^2/2)_x + (u^2/2)_y = 0$, subject to the initial data
 $u(x,t\!=\!0) = \sin^2(\pi x) \sin^2(\pi y)$ and periodic boundary
conditions on $[0,1]^2$.  In Figure \ref{figure:2d_burgers} we present
the solution obtained at time $T=1.5$ with different mesh spacings.
One can easily notice that the shocks are well resolved and there
are no spurious oscillations.

\begin{figure}[htbp]
  \begin{center}
  \includegraphics[width=4.5in]{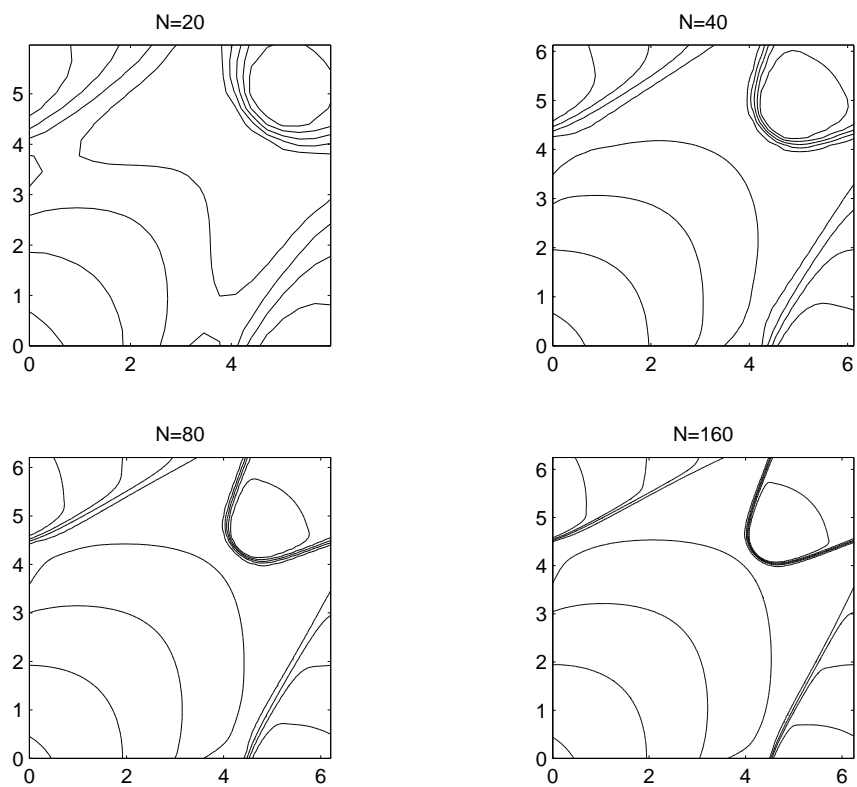}
  \end{center}
  \caption{Two-Dimensional Burgers equation, $\lambda=0.425$, $T=1.5$
 \label{figure:2d_burgers}}
\end{figure}

\end{enumerate}


\end{document}